\documentclass[a4paper,11pt] {amsart}
\usepackage{amsthm,amssymb,latexsym,mathrsfs}
\input amssym.def

\author{Beno\^it F. Sehba}

\title[Maximal function and Carleson measures ]{Maximal functions and measures on the upper-half plane}
\newtheorem{theorem}{T{\hskip 0pt\footnotesize\bf HEOREM}}[section]
\newtheorem{lemma}[theorem]{L{\hskip 0pt\footnotesize\bf EMMA}}
\newtheorem{proposition}[theorem]{P{\hskip 0pt\footnotesize\bf ROPOSITION}}
\newtheorem{definition}[theorem]{D{\hskip 0pt\footnotesize\bf EFINITION}}
\newtheorem{corollary}[theorem]{C{\hskip 0pt\footnotesize\bf OROLLARY}}

\newcommand{\bprop} {\begin{proposition}}
\newcommand{\eprop} {\end{proposition}}
\newcommand{\btheo} {\begin{theorem}}
\newcommand{\etheo} {\end{theorem}}
\newcommand{\blem} {\begin{lemma}}
\newcommand{\elem} {\end{lemma}}
\newcommand{\bcor} {\begin{corollary}}
\newcommand{\ecor} {\end{corollary}}

\newcommand{\Be}{\begin{equation}}
\newcommand{\Ee}{\end{equation}}
\newcommand{\Bea}{\begin{eqnarray}}
\newcommand{\Eea}{\end{eqnarray}}
\newcommand{\Bes}{\begin{equation*}}
\newcommand{\Ees}{\end{equation*}}
\newcommand{\Beas}{\begin{eqnarray*}}
\newcommand{\Eeas}{\end{eqnarray*}}
\newcommand{\Ba}{\begin{array}}
\newcommand{\Ea}{\end{array}}

\scrollmode
\begin{document}
%\date{\today}
%\address{Carnot D. Kenfack, D\'epartement de math\'ematiques, Facult\'e des Sciences, Universit\'e de Yaound\'e I, BP.812 Yaound\'e, Cameroun}
%\email{dondjiocarnot@yahoo.fr}
%\address{\textit{Recent address}: Gonessa Jocelyn, LATP, UM.R. CNRS. 6632, CMI, Universit\'e de Provence, 39 Rue F-Juliot-Curie 13453, Marseille Cedex 13, %France}
%\email{gonessa@cmi.univ-mrs.fr}
\address{Beno\^it F. Sehba, Department of Mathematics, University of Ghana, Legon, P. O. Box LG 62, Legon, Accra, Ghana}
\email{bfsehba@ug.edu.gh}
%\address{Jocelyn Gonessa, .}
%\email{j_gonessa@yahoo.fr}
%\address{Edgar Tchoundja, Centre de Recerca Matem\`atica, 08193 Bellaterra Barcelona, Spain.}
%\email{etchoundja@crm.cat}
%\keywords{Hankel operators, Hardy-Orlicz spaces,  weighted BMOA spaces.}
%\subjclass[2000]{Primary 47B35, Secondary 32A35, 32A37}

%\date{today}
%\address{Jocelyn Gonessa, .}
%\email{j_gonessa@yahoo.fr}
%\address{Beno\^it Florent Sehba, Centre de Math\'ematiques et Informatique (CMI),
%Universit\'e de Provence, Technop\^ole Ch\^ateau-Gombert,
%39, rue F. Joliot Curie,
%13453 Marseille Cedex 13
%FRANCE .}
%\email{sehbab@cmi.univ-mrs.fr}
\keywords{B\'ekoll\`e-Bonami weight, Carleson-type embedding, Dyadic grid, Maximal function, Upper-half plane.}
\subjclass[2000]{Primary: 47B38 Secondary: 30H20, 42A61, 42C40}
%\thanks{The second author was partially supported by the ANR project ANR-09-BLAN-0058-01}

%\noindent  ABSTRACT.\,

\maketitle

\begin{abstract}
We study weighted boundedness of Hardy-Littlewood-type maximal function involving Orlicz functions. We also obtain some sufficient conditions for the weighted boundedness of the Hardy-Littlewood maximal function of the upper-half plane.
\end{abstract}

%\part{Use this type of header for very long papers only}
% use lowercase except for proper names

\section{Introduction}
Our setting is the upper-half plane $$\mathcal{H}=\{z=x+iy\in \mathbb {C}: x\in \mathbb R,\,\,\, \textrm{and}\,\,\,y>0\}.$$ For $\omega$ a weight, that is a nonnegative locally integrable function on $\mathcal H$, $\alpha>-1$, and $1\le p<\infty$, we denote by $L_\omega^p(\mathcal H, dV_\alpha)$, the set of functions $f$ defined on $\mathcal H$ such that
$$\|f\|_{p,\omega,\alpha }^p:=\int_{\mathbb D}|f(z)|^p\omega(z)dV_\alpha(z)<\infty$$
with $dV_\alpha(x+iy)=y^\alpha dxdy$. We write $L^p(\mathcal H, dV_\alpha)$ when $\omega(z)=1$ for any $z\in \mathcal{H}$, and put $\|f\|_{p,\alpha }=\|f\|_{p,1,\alpha }$. We recall that the Bergman space $A_\alpha^p(\mathcal{H})$ is the closed subset of $L^p(\mathcal H, dV_\alpha)$ consisting of holomorphic functions.
\vskip .2cm
Let us recall that for any interval $I\subset \mathbb{R}$, its associated Carleson square $Q_I$ is the set $Q_I:=\{z=x+iy\in \mathbb C:x\in I\,\,\,\textrm{and}\,\,\,0<y<|I|\}$.
%\vskip .1cm
%Given a function $\Phi:[0,\infty)\rightarrow [0,\infty)$, We say  $\Phi$ is a growth function if it is a continuous and non-decreasing function.
\vskip .1cm
Let $\Phi$ be a Young function, and $\alpha>-1$. For any interval $I\subset \mathbb{R}$, define $L^\Phi(Q_I, |Q_I|_{\omega,\alpha}^{-1}\omega dV_\alpha)$ to be the space of all functions $f$ such that 
$$\frac{1}{|Q_I|_{\omega,\alpha}}\int_{Q_I}\Phi\left(|f(z)|\right)\omega(z)dV_\alpha(z)<\infty$$
where $|Q_I|_{\omega,\alpha}=\int_{Q_I}\omega(z)dV_\alpha(z)$.
We define on $L^\Phi(Q_I, |Q_I|_{\omega,\alpha}^{-1}\omega dV_\alpha)$ the following Luxembourg norm $$\|f\|_{Q_I,\Phi,\omega,\alpha}:=\inf\{\lambda>0: \frac{1}{|Q_I|_{\omega,\alpha}}\int_{Q_I}\Phi\left(\frac{|f(z)|}{\lambda}\right)\omega(z)dV_\alpha(z)\le 1\}.$$
 When $\Phi(t)=t^p$, $1\le p<\infty$, $L^\Phi(Q_I, |Q_I|_{\omega,\alpha}^{-1}\omega dV_\alpha)$ is just $L^p(Q_I, |Q_I|_{\omega,\alpha}^{-1}\omega dV_\alpha)$ in which case $\|f\|_{Q_I,\Phi,\omega,\alpha}$ is just replaced by $\left(\frac{1}{|Q_I|_{\omega,\alpha}}\int_{Q_I}|f(z)|^p\omega(z)dV_\alpha(z)\right)^{1/p}$. Then the maximal function $\mathcal{M}_{\Phi,\omega,\alpha}$ is defined as $$\mathcal{M}_{\Phi,\omega,\alpha}f(z):=\sup_{z\in Q_I}\|f\|_{Q_I,\Phi,\omega,\alpha}.$$
More precisely, the supremum is taken over all intervals $I$ such that $z\in Q_I$ . Our definition here is inspired from the one in the works \cite{Bagby,Perez} and actually, it is a weighted version of the one studied in \cite{Perez}. When $\omega=1$, we simply write $\|f\|_{Q_I,\Phi,\alpha}$ for $\|f\|_{Q_I,\Phi,1,\alpha}$ and $\mathcal{M}_{\Phi,\alpha}$ for the corresponding maximal function.
We observe that when $\Phi=1$, $\mathcal{M}_{\Phi,\omega,\alpha}$ coincides with the (weighted) Hardy-Littlewood maximal function $\mathcal{M}_{\omega,\alpha}$ given by $$\mathcal{M}_{\omega,\alpha} f(z)=\sup_{I\subset \mathbb{R}}\frac{\chi_{Q_I}(z)}{|Q_I|_{\omega,\alpha}}\int_{Q_I}|f(z)|\omega(z)dV_\alpha(z).$$
The unweighted Hardy-Littlewood maximal function, that corresponds to $\omega=1=\Phi$ will be denoted $\mathcal{M}_{\alpha}$.
\vskip .2cm
We are interested in this paper in weighted boundedness of the maximal function $\mathcal{M}_{\Phi,\omega,\alpha}$, that is the characterization of positive Borel measure $\mu$ and weight $\omega$ such that 
$\mathcal{M}_{\Phi,\omega,\alpha}$ is bounded from $L^p(\mathcal H, \omega dV_\alpha)$  to $L^q(\mathcal H, d\mu)$. These estimates are quite useful in the estimate of other operators as the Bergman projection (see for example \cite{AlPottReg,PR}). We will also see that they can be used to obtain alternative characterizations of Bergman spaces. We observe that in \cite{CarnotBenoit}, these characterizations were obtained for the weighted Hardy-Littlewood maximal function of the upper-half plane. 
\subsection{Statement of the results}
Recall that a function from $[0,\infty)$ to itself is a Young function if it is continuous, convex and increasing, and satisfies $\Phi(0)=0$ and $\Phi(t)\rightarrow \infty$ as $t\rightarrow \infty$.  We assume all over the text that the Young function $\Phi$ is such that $\Phi(1)=1$. We also recall that the function $t\mapsto \frac{\Phi(t)}{t}$ is increasing, and $\Phi'(t)\simeq \frac{\Phi(t)}{t}$.  
\vskip .3cm
Given a Young function $\Phi$, we say it satisfies the $\Delta_2$ (or doubling) condition, if there exists a constant $K>1$ such that, for any $t\ge 0$,
\begin{equation}\label{eq:delta2condition}
 \Phi(2t)\le K\Phi(t).
 \end{equation}

Let $1<p<\infty$. We say a Young function $\Phi$ belongs to the class $B_p$, if it satisfies the $\Delta_2$ condition and there is a positive constant $c$ such that 
\Be\label{eq:BpDini}
\int_c^\infty \frac{\Phi(t)}{t^p}\frac{dt}{t}<\infty.
\Ee
For any set $E\subset \mathcal{H}$, given a weight $\omega$, and $\alpha>-1$, we write $$|E|_{\omega,\alpha}:=\int_E\omega(z)dV_\alpha(z).$$
%Our first result is as follows.

Recall that for $1<p<\infty$ and $\alpha>-1$, a weight $\omega$ is said to belong to the B\'ekoll\`e-Bonami class $\mathcal{B}_{p,\alpha}$,  if $[\omega]_{\mathcal{B}_{p,\alpha}}<\infty$, where
$$[\omega]_{\mathcal{B}_{p,\alpha}}:=\sup_{I\subset \mathbb R,\,\,\, I\,\,\,\textrm{interval} }\left(\frac{1}{|Q_I|_\alpha}\int_{Q_I}\omega(z)dV_\alpha(z)\right)\left(\frac{1}{|Q_I|_\alpha}\int_{Q_I}\omega(z)^{1-p'}dV_\alpha(z)\right)^{p-1}.$$
For $p=1$, a weight $\omega$ is said to belong to the B\'ekoll\`e-Bonami class $\mathcal{B}_{1,\alpha}$,  if $$[\omega]_{\mathcal{B}_{1,\alpha}}:=[\omega]_{\mathcal{B}_{1,\alpha}}:=\sup_I\mbox{ess}\sup_{z\in Q_I}\left(\frac{1}{|Q_I|_{\alpha}}\int_{Q_I}\omega dV_{\alpha}\right)\omega(z)^{-1}<\infty.$$
It is easy to see that $\mathcal{B}_{1,\alpha}$ is a subset of $\mathcal{B}_{2,\alpha}$ and $$[\omega]_{\mathcal{B}_{2,\alpha}}\le [\omega]_{\mathcal{B}_{1,\alpha}},\,\,\,\textrm{for all}\,\,\,\omega\in \mathcal{B}_{1,\alpha}.$$
We define $$\mathbb{B}_{\infty,\alpha}:=\bigcup_{p>1}\mathcal{B}_{p,\alpha}.$$
%For $\alpha>-1$, we say a weight $\omega$ is $\alpha$-doubling, if there are an increasing function $\varphi$ with $\varphi(1)=1$ and a constant $K=K(\omega,\alpha)\geq 1$ such that for any interval $I\subset \mathbb{R}$ and any set $E\subset Q_I$, \Be\label{eq:alphadoubling}\frac{|Q_I|_{\omega,\alpha}}{|E|_{\omega,\alpha}}\le K\varphi\left(\frac{|Q_I|_\alpha}{|E|_\alpha}\right).\Ee
 
We have the following result for the weighted maximal function.
\btheo\label{thm:main1}
Let $\alpha>-1$, $1<p\le q<\infty$, and $\omega$ a weight and $\mu$ a positive Borel measure on $\mathcal H$. Assume that $\Phi\,\in B_{p}$ and that $\omega\in \mathbb{B}_{\infty,\alpha}$. Then the following assertions are equivalent.
 \begin{itemize}
\item[(i)] There exists a constant $C_1>0$ such that for any $f\in L^p(\mathcal H, \omega dV_\alpha)$,
\Be\label{eq:main11}
\left(\int_{\mathcal H}\left(\mathcal{M}_{\Phi,\omega,\alpha} f(z)\right)^qd\mu(z)\right)^{1/q}\le C_1\left(\int_{\mathcal H}|f(z)|^p\omega(z)dV_\alpha(z)\right)^{1/p}.
\Ee

\item[(ii)] There is a constant $C_2$ such that for any interval $I\subset \mathbb R$,
\Be\label{eq:main12} \mu(Q_I)\leq C_2 |Q_I|_{\omega,\alpha}^{\frac{q}{p}}.\Ee
%$|Q_I|_{\omega,\alpha}=\int_{Q_I}\omega dV_\alpha$.
\end{itemize}
\etheo
We have the following lower triangle case result.
\btheo\label{thm:main2}
Let $\alpha>-1$, $1<q<p<\infty$. Let $\omega$ be a weight and $\mu$ a positive Borel measure on $\mathcal H$. Assume that $\Phi\,\in B_{p}$ and that $\omega\in \mathbb{B}_{\infty,\alpha}$. Then (\ref{eq:main11}) holds if and only if the function
\Be\label{eq:main21} K_{\mu}(z):=\sup_{I\subset \mathbb R,\,\,\,I\,\,\,\textrm{interval},\,\,\,z\in Q_I}\frac{\mu(Q_{I})}{|Q_I|_{\omega,\alpha}}
\Ee belongs to
$L^s(\mathcal H,\omega dV_\alpha)$ where $s=\frac{p}{p-q}$.
\etheo
We recall that the complementary function $\Psi$ of the Young function $\Phi$, is the function defined from $\mathbb R_+$ onto itself by
\begin{equation}\label{complementarydefinition}
\Psi(s)=\sup_{t\in\mathbb R_+}\{ts - \Phi(t)\}.
\end{equation}
We remark that a Young function satisfies (\ref{eq:BpDini}) if and only if its complementary function $\Psi$ satisfies
\Be\label{eq:BpDiniconj}
\int_c^\infty \left(\frac{t^{p'}}{\Psi(t)}\right)^{p-1}\frac{dt}{t}<\infty.
\Ee
Here, and all over the text, $p'$ denotes the conjugate exponent of $p$. We recall that 
\Be\label{phipsiidentity}
t\le \Phi^{-1}(t)\Psi^{-1}(t)\le 2t,\,\,\,\textrm{for all}\,\,\,t>0.
\Ee
The following result provides a sufficient condition for the off-diagonal boundedness of the Hardy-Littlewood maximal function $\mathcal{M}_\alpha$.
\btheo\label{thm:main3}
Let $\alpha>-1$, $1<p\le q<\infty$. Let $\Phi\in B_p$ and denote by $\Psi$ its complementary function. Assume that $\omega$ is a weight and $\mu$ a positive Borel measure on $\mathcal H$ such that there is positive constant $C$ for which for any interval $I\subset \mathbb{R}$,
\Be\label{eq:main31}
\|\omega^{-1}\|_{Q_I,\Psi,\alpha}^q\mu(Q_I)\le C|Q_I|_\alpha^{q/p}.
\Ee

 Then there is a positive constant $K$ such that for any $f\in L^p(\mathcal H, \omega dV_\alpha)$, 
 \Be\label{eq:main32}
 \left(\int_{\mathcal H}\left(\mathcal{M}_{\alpha} f(z)\right)^qd\mu(z)\right)^{1/q}\le K\|f\omega\|_{p,\alpha}.
 \Ee
\etheo
It is easy to see that for $1<p<\infty$, and $r>1$, $\Phi(t)=t^{(p'r)'}$ is in the class $B_p$. Thus we have the following corollary.
\bcor\label{cor:main1}
Let $\alpha>-1$, $1<p\le q<\infty$.  Assume that $\omega$ is a weight and $\mu$ a positive Borel measure on $\mathcal H$ such that for some $r>1$, there is positive constant $C$ for which for any interval $I\subset \mathbb{R}$,
\Be\label{eq:main310}
\left(\frac{1}{|Q_I|_\alpha}\int_{Q_I}\omega^{-p'r}\right)^{q/p'r}\mu(Q_I)\le C|Q_I|_\alpha^{q/p}.
\Ee

 Then there is a positive constant $K$ such that for any $f\in L^p(\mathcal H, \omega dV_\alpha)$, 
 \Be\label{eq:main320}
 \left(\int_{\mathcal H}\left(\mathcal{M}_{\alpha} f(z)\right)^qd\mu(z)\right)^{1/q}\le K\|f\omega\|_{p,\alpha}.
 \Ee
\ecor
We also obtain the following result.
\btheo\label{thm:main4}
Let $\alpha>-1$, $1<p\le q<\infty$. Let $\mu$ be a positive measure and $\omega$ a weight such that 
there is a positive constant $C_1$ such that for any interval $I\subset \mathbb R$,
$$ \mu(Q_I)\leq C_1 |Q_I|_{\omega,\alpha}^{\frac{q}{p}}.$$ Then there is a constant $C_2>0$ such that for any function $f$,

\Be\label{eq:main41}
\left(\int_{\mathcal H}\left(\mathcal{M}_{\Phi,\alpha} f(z)\right)^q d\mu(z)\right)^{1/q}\le C\left(\int_{\mathcal H}|f(z)|^p\mathcal{M}_\alpha \omega(z)d\mu(z)\right)^{1/p}
\Ee
\etheo
%An easy consequence of the above result is that an analytic function $f$ belongs to the Bergman space $A_\omega^p(\mathcal H)$ if and only if $M_\omega f\in L_\omega^p(\mathcal H)$. 
We also have the following.

\btheo\label{thm:main5}
Let $\alpha>-1$, $1<q<p<\infty$. Let $\omega$ be a weight and $\mu$ a positive Borel measure on $\mathcal H$. Assume that $\Phi\,\in B_{p}$ and that  the function
\Be\label{eq:main51} K_{\mu}(z):=\sup_{I\subset \mathbb R,\,\,\,I\,\,\,\textrm{interval},\,\,\,z\in Q_I}\frac{\mu(Q_{I})}{|Q_I|_{\omega,\alpha}}
\Ee belongs to
$L^s(\mathcal H,(\mathcal{M}_\alpha\omega) dV_\alpha)$ where $s=\frac{p}{p-q}$. Then (\ref{eq:main41}) holds for any function $f$.
\etheo
\subsection{Methods of proof}
Our presentation is essentially based on discretrization methods with some of our considerations been different from the ones in \cite{Perez} where similar topics were considered in $\mathbb{R}^n$. For the proof of Theorem \ref{thm:main1}, we will prove the following inclusion, $$\{z\in \mathcal{H}: \mathcal{M}_{\Phi,\omega,\alpha}f(z)>\lambda\}\subset \{z\in \mathcal{H}: \mathcal{M}_{\Phi,\omega,\alpha}^df(z)>\frac{\lambda}{C_\alpha}\}$$ 
where $\mathcal{M}_{\Phi,\omega,\alpha}^d$ is the dyadic analogue of $\mathcal{M}_{\Phi,\omega,\alpha}$. We also use the usual covering of  $\{z\in \mathcal{H}: \mathcal{M}_{\Phi,\omega,\alpha}^df(z)>\frac{\lambda}{C_\alpha}\}$ and prove that under condition (\ref{eq:main12}), $$\mu\left(\{z\in \mathcal{H}: \mathcal{M}_{\Phi,\alpha}^df(z)>\frac{\lambda}{C_\alpha}\}\right)\le C|\{z\in \mathcal{H}: \mathcal{M}_{\Phi,\alpha}^df(z)>\frac{\lambda}{C_\alpha}\}|_{\omega,\alpha}^{q/p}.$$
\vskip .3cm
In the proof of Theorem \ref{thm:main2}, we use a discretization of the integral and an extension of the Carleson Embedding Theorem. The Carleson Embedding Theorem needed here is a generalization of the classical one (see \cite{Pereyra,Sehba}) as we replace the average $\frac{1}{|Q_I|_{\omega,\alpha}}\int_{Q_I}|f(z)|\omega(z)dV_\alpha(z)$ by $\|f\|_{Q_I,\Phi,\omega,\alpha}$.
\vskip .3cm
For the proof of Theorem \ref{thm:main4} and Theorem \ref{thm:main5}, we also use discretization and the unweighted version of Theorem \ref{thm:main1}.
\section{Proof of Theorem \ref{thm:main1}}
For $\alpha>-1$, we say a weight $\omega$ is $\alpha$-doubling, if there are an increasing function $\varphi$ with $\varphi(1)=1$ and a constant $K=K(\omega,\alpha)\geq 1$ such that for any interval $I\subset \mathbb{R}$ and any set $E\subset Q_I$, \Be\label{eq:alphadoubling}\frac{|Q_I|_{\omega,\alpha}}{|E|_{\omega,\alpha}}\le K\varphi\left(\frac{|Q_I|_\alpha}{|E|_\alpha}\right).\Ee

Let us start this section with the following lemma.
\blem\label{lem:covering}
Let $\alpha>-1$, let $\omega$ be a weight on $\mathcal{H}$, and assume that $\Phi$ is a Young function. Then for any compactly supported function $f$ and any $\lambda>0$, there exists a family of maximal (with respect to inclusion) dyadic intervals $\{I_j\}_j$ such that $$\{z\in \mathcal{H}: \mathcal{M}_{\Phi,\omega,\alpha}^df(z)>\lambda\}=\bigcup_jQ_{I_j}.$$ The above Carleson squares $Q_{I_j}$ are maximal such $$\|f\|_{Q_{I_j},\Phi,\omega,\alpha}>\lambda\,\,\,\textrm{for each}\,\,\,j. $$
If $\omega$ is $\alpha$-doubling, then $$\|f\|_{Q_{I_j},\Phi,\omega,\alpha}\le \varphi(2^{2+\alpha}) K\lambda\,\,\,\textrm{for each}\,\,\,j$$ where $\varphi$ and $K$ are the function and the constant in (\ref{eq:alphadoubling}).
Moreover, for any weight $\omega$,
$$\left|\{z\in \mathcal{H}:\mathcal{M}_{\Phi,\omega,\alpha}^df(z)>\lambda\}\right|_{\omega,\alpha}\le C\int_{\{z\in \mathcal{H}:|f(z)|>\frac{\lambda}{2}\}}\Phi\left(\frac{|f(z)|}{\lambda}\right)\omega(z)dV_\alpha(z).$$
\elem
\begin{proof}
That $\{z\in \mathcal{H}: \mathcal{M}_{\Phi,\omega,\alpha}^df(z)>\lambda\}$ is a union of maximal dyadic Carleson squares $Q_{I_j}$ such that $\|f\|_{Q_{I_j},\Phi,\omega,\alpha}>\lambda$ follows from the usual arguments. 
\vskip .3cm
Now assume that $\omega$ is $\alpha$-doubling.
To see that $\|f\|_{Q_{I_j},\Phi,\omega,\alpha}\le K\varphi(2^{2+\alpha})\lambda$, observe that if $I$ is such that $Q_I$ is one of the maximal Carleson squares above, then $\|f\|_{Q_{\tilde{I}},\Phi,\omega,\alpha}\le \lambda$, where $\tilde{I}$ is a parent of $I$. Using the convexity of $\Phi$, we obtain that for any $t>0$,
\Beas
L &:=& \frac{1}{|Q_I|_{\omega,\alpha}}\int_{Q_I}\Phi\left(\frac{|f(z)|}{\varphi(2^{2+\alpha})Kt}\right)\omega(z)dV_\alpha(z)\\ &\le&
\frac{\varphi(2^{2+\alpha})K}{|Q_{\tilde{I}}|_{\omega,\alpha}}\int_{Q_{\tilde{I}}}\Phi\left(\frac{|f(z)|}{\varphi(2^{2+\alpha})Kt}\right)\omega(\omega)dV_\alpha(z)\\ &\le& \frac{1}{|Q_{\tilde{I}}|_{\omega,\alpha}}\int_{Q_{\tilde{I}}}\Phi\left(\frac{|f(z)|}{t}\right)\omega(z)dV_\alpha(z).
\Eeas
Thus $\|f\|_{Q_{I},\Phi,\omega,\alpha}\le \varphi(2^{2+\alpha})K\|f\|_{Q_{\tilde{I}},\Phi,\omega,\alpha}$.
\vskip .3cm
Note that as for any Carleson square $Q_{I_j}$ in the above family, we have $\|f\|_{Q_{I_j},\Phi,\omega,\alpha}>\lambda$, it comes that $$\frac{1}{|Q_{I_j}|_{\omega,\alpha}}\int_{Q_{I_j}}\Phi\left(\frac{|f(z)|}{\lambda}\right)\omega(z)dV_\alpha(z)>1\,\,\,\textrm{for all}\,\,\,j.$$ That is $$|Q_{I_j}|_{\omega,\alpha}<\int_{Q_{I_j}}\Phi\left(\frac{|f(z)|}{\lambda}\right)\omega(z)dV_\alpha(z)\,\,\,\textrm{for all}\,\,\,j.$$
Hence
\Beas
\left|\{z\in \mathcal{H}:\mathcal{M}_{\Phi,\omega,\alpha}^df(z)>\lambda\}\right|_{\omega,\alpha} &=& \left|\bigcup_{j}Q_{I_j}\right|_{\omega,\alpha}=\sum_{j}|Q_{I_j}|_{\omega,\alpha}\\
&\le& \sum_{j}\int_{Q_{I_j}}\Phi\left(\frac{|f(z)|}{\lambda}\right)\omega(z)dV_\alpha(z)\\ &\le& \int_{\mathcal{H}}\Phi\left(\frac{|f(z)|}{\lambda}\right)\omega(z)dV_\alpha(z).
\Eeas
Now following the usual arguments, we write $f=f_1+f_2$ where $$f_1:=f\chi_{\{z\in \mathcal{H}:|f(z)|>\frac{\lambda}{2}\}}.$$ Then 
$$\mathcal{M}_{\Phi,\omega,\alpha}^df(z)\le \mathcal{M}_{\Phi,\omega,\alpha}^df_1(z)+\mathcal{M}_{\Phi,\omega,\alpha}^df_2(z)\le \mathcal{M}_{\Phi,\omega,\alpha}^df_1(z)+\frac{\lambda}{2}.$$
It follows that

\Beas\left|\{z\in \mathcal{H}:\mathcal{M}_{\Phi,\omega,\alpha}^df(z)>\lambda\}\right|_{\omega,\alpha} &\le& \left|\{z\in \mathcal{H}:\mathcal{M}_{\Phi,\omega,\alpha}^df_1(z)>\frac{\lambda}{2}\}\right|_{\omega,\alpha}\\ &\le& C\int_{\{z\in \mathcal{H}:|f(z)|>\frac{\lambda}{2}\}}\Phi\left(\frac{|f(z)|}{\lambda}\right)\omega(z)dV_\alpha(z).
\Eeas
The proof is complete.
\end{proof}

Let us prove the following level sets embedding.

\blem\label{lem:levelsetsembed}
Let $\alpha>-1$ and let $\Phi$ be a Young function. Assume that $\omega$ is an $\alpha$-doubling weight. Let $f$ be a locally integrable function. Then for any $\lambda>0$, 
\Be\label{eq:levelsetsembed}
\{z\in \mathcal H: \mathcal{M}_{\Phi,\omega,\alpha}f(z)>\lambda\}\subset \{z\in \mathcal{H} \mathcal{M}_{\Phi,\omega,\alpha}^df(z)>\frac{\lambda}{C}\}
\Ee
with $C=2K\varphi(2^{2+\alpha})(1+(K\varphi(2^{2+\alpha}))^2)$, where $\varphi$ and $K$ are the function and the constant in (\ref{eq:alphadoubling}).
\elem
\begin{proof}
The proof is essentially the same as the one of Lemma 3.4 in \cite{CarnotBenoit}. We give here the main modifications.
Let us put $$E_\lambda :=\{z\in \mathcal H: \mathcal{M}_{\Phi,\omega,\alpha}f(z)>\lambda\}$$ and $$E_{\lambda}^d:=\{z\in \mathcal H: \mathcal{M}_{\Phi,\omega,\alpha}^df(z)>\frac{\lambda}{C}\}.$$ Recall that there is a family $\{Q_{I_j}\}_{j\in \mathbb N_0}$ of maximal (with respect to the inclusion) dyadic Carleson squares such that
$$\frac{\lambda}{C}<\|f\|_{Q_{I_j},\Phi,\alpha}\le K\varphi(2^{2+\alpha})\frac{\lambda}{C}\,\,\,\textrm{for each}\,\,\,j $$
$K$ being the constant in (\ref{eq:alphadoubling}). Moreover, $E_{\lambda}^d=\cup_{j\in \mathbb N}Q_{I_j}$.
\vskip .1cm
Let $z\in E_{\lambda}$ and suppose that $z\notin E_{\lambda}^d$. Note   that there is an interval $I$ (not necessarily dyadic) such that $z\in Q_I$ and
\Be\label{eq:hypt}\|f\|_{Q_{I},\Phi,\omega,\alpha} >\lambda.\Ee
Recall with \cite[Lemma 2.3.]{CarnotBenoit} that $I$ can be covered by at most two adjacent dyadic intervals $J_1$ and $J_2$ ($J_1$ on the left of $J_2$) such that $|I|<|J_1|=|J_2|\le 2|I|$. Hence  $Q_I\subset Q_{J_1}\cup Q_{J_2}$ and we have that $z$ belongs only to one and only one of the Carleson boxes $Q_{J_1}$ and $Q_{J_2}$. Let us suppose that $z\in Q_{J_1}$ (in which case $z\notin Q_{J_2}$). Then  $Q_{J_1}\cap E_{\lambda}^d =\emptyset$ or $Q_{J_1}\supset Q_{I_j}$ for some $j$ and in both cases, because of the maximality of the $I_j$s, we obtain that
$$\|f\|_{Q_{J_1},\Phi,\omega,\alpha}\le \frac{\lambda}{C}.$$
For the other interval $J_2$, we have the following possibilities
%$$
%\left\{ \begin{matrix} J_2=I_j\quad \textrm{for some}\quad j\\
%
%               J_2\subset I_j \quad \textrm{for some}\quad j\\
%               J_2\supset I_j\quad \textrm{for some}\quad j\\
%               J_2\cap B=\emptyset.
%                              \end{matrix}\right.$$
\begin{eqnarray*}
&~& J_2=I_j\quad \textrm{for some}\quad j;\\
&~& J_2\subset I_j \quad \textrm{for some}\quad j;\\
&~& J_2\supset I_j\quad \textrm{for some}\quad j;\\
&~& J_2\cap B=\emptyset.
\end{eqnarray*}

If $J_2\supset I_j$ for some $j$ or $J_2\cap E_{\lambda}^d=\emptyset$, then because of the maximality of the $I_j$s,

\[\|f\|_{Q_{J_2},\Phi,\omega,\alpha}\le \frac{\lambda}{C}.\]

If $J_2=I_j$ for some $j$, then 
$$\|f\|_{Q_{J_2},\Phi,\omega,\alpha}\le \frac{K\lambda\varphi(2^{2+\alpha})}{C}.$$
 If $J_2\subset I_j$ for some $j$, then 
%\[
%\left\{ \begin{matrix} J_2=I_j^-\\ J_2\subset I_j^- \\  J_2\subseteq I_j^+ \end{matrix}\right.\]
\[J_2=I_j^-,\,\,\, J_2\subset I_j^-\,\,\,\textrm{or}\,\,\,   J_2\subseteq I_j^+\]
where $I_j^-$ and $I_j^+$ denote the left and right halfs of $I_j$ respectively. If $J_2\subset I_j^-$ or $J_2\subseteq I_j^+$, then $J_1\cap I_j\neq \emptyset$, and this implies that $J_1\subset I_j$. Thus $z\in Q_{J_1}\subset Q_{I_j}\subset E_{\lambda}^d$ which contradicts the hypothesis $z\notin E_{\lambda}^d$. Hence we can only have $J_2=I_j^-$ which leads to the estimate
$$\|f\|_{Q_{J_2},\Phi,\omega,\alpha}\le K\varphi(2^{2+\alpha})\|f\|_{Q_{I_j},\Phi,\omega,\alpha}\le \frac{(K\varphi(2^{2+\alpha}))^2}{C}\lambda.$$

It follows from that above discussion and the convexity of $\Phi$ that
\Beas
L &:=& \frac{1}{|Q_{I}|_{\omega,\alpha}}\int_{Q_{I}}\Phi\left(\frac{|f|}{2K\varphi(2^{2+\alpha})\left(\frac{1}{C}+\frac{(K\varphi(2^{2+\alpha}))^2}{C}\right)\lambda}\right)\omega dV_\alpha \\ &\le&
\frac{|Q_{J_1}|_{\omega,\alpha}}{|Q_{I}|_{\omega,\alpha}}\left(\frac{1}{|Q_{J_1}|_{\omega,\alpha}}\int_{Q_{J_1}}\Phi\left(\frac{|f|}{2K\varphi(2^{2+\alpha})\frac{\lambda}{C}}\right)\omega dV_\alpha\right)+\\ &&\frac{|Q_{J_2}|_{\omega,\alpha}}{|Q_{I}|_{\omega,\alpha}}\left(\frac{1}{|Q_{J_2}|_{\omega,\alpha}}\int_{Q_{J_2}}\Phi\left(\frac{|f|}{2(K\varphi(2^{2+\alpha}))^3\frac{\lambda}{C}}\right)\omega dV_\alpha\right)\\  &\le& \frac{K\varphi(2^{2+\alpha})}{|Q_{J_1}|_{\omega,\alpha}}\int_{Q_{J_1}}\Phi\left(\frac{|f|}{2K\varphi(2^{2+\alpha})\frac{\lambda}{C}}\right)\omega dV_\alpha+\\ && \frac{K\varphi(2^{2+\alpha})}{|Q_{J_2}|_{\omega,\alpha}}\int_{Q_{J_2}}\Phi\left(\frac{|f|}{2(K\varphi(2^{2+\alpha}))^3\frac{\lambda}{C}}\right)\omega dV_\alpha\\  
&\le& \frac{1}{2}\left(\frac{1}{|Q_{J_1}|_{\omega,\alpha}}\int_{Q_{J_1}}\Phi\left(\frac{|f|}{\frac{\lambda}{C}}\right)\omega dV_\alpha\right)+\\ && \frac{1}{2}\left(\frac{1}{|Q_{J_2}|_{\omega,\alpha}}\int_{Q_{J_2}}\Phi\left(\frac{|f|}{(K\varphi(2))^2\frac{\lambda}{C}}\right)\omega dV_\alpha\right)\\ &\le& 1.
\Eeas
Thus 
$$\|f\|_{Q_{I},\Phi,\omega,\alpha}\le 2\lambda\left(\frac{K\varphi(2^{2+\alpha})}{C}+\frac{(K\varphi(2^{2+\alpha}))^3}{C}\right)=\lambda.$$
which clearly contradicts (\ref{eq:hypt}). The proof is complete.
%\end{eqnarray*}
\end{proof}
We need the following estimate.
\blem\label{lem:carlmeasineq}
Let $\gamma \geq 1$, and $\alpha>-1$. Let $\sigma$ and $\omega$ be weights, and $\mu$ a positive measure on $\mathcal{H}$. Assume that there is a constant $C>0$ such that for any interval $I\subset \mathbb{R}$,
$$\mu(Q_I)\le C|Q_I|_{\omega,\alpha}^\gamma.$$
Then  for any function $f$ and any $t>0$, $$\mu\left(\{z\in \mathcal{H}:\mathcal{M}_{\Phi,\sigma,\alpha}^df(z)>t\}\right)\le C\left|\left(\{z\in \mathcal{H}:\mathcal{M}_{\Phi,\sigma,\alpha}^df(z)>t\}\right)\right|_{\omega,\alpha}^\gamma.$$
\elem
\begin{proof}
Recall with Lemma \ref{lem:covering} that there is family $\{I_j\}_j$ of dyadic maximal intervals such that $$\{z\in \mathcal{H}:\mathcal{M}_{\Phi,\sigma,\alpha}^df(z)>t\}=\bigcup_j Q_{I_j}.$$ Thus 
\Beas
\mu\left(\{z\in \mathcal{H}:\mathcal{M}_{\Phi,\sigma,\alpha}^df(z)>t\}\right) &=& \sum_j\mu(Q_{I_j})\\
&\le& C\sum_j|Q_{I_j}|_{\omega,\alpha}^\gamma\\ &\le& C\left(\sum_j|Q_{I_j}|_{\omega,\alpha}\right)^\gamma \\ &=& C\left|\{z\in \mathcal{H}:\mathcal{M}_{\Phi,\sigma,\alpha}^df(z)>t\}\right|_{\omega,\alpha}^\gamma 
\Eeas
\end{proof}
\subsection{A first proof of Theorem \ref{thm:main1}}
The following result extends the classical estimate of the weighted Hardy-Littlewood maximal function.
\blem\label{lem:estimmaxphi}
Let $1<p<\infty$ and $\alpha>-1$. Assume that $\omega$ is a weight and $\Phi\in B_p$. Then there is a positive constant $C$ such that for any function $f$,
\Be\label{eq:estmaxphi}
\left(\int_{\mathcal{H}}(\mathcal{M}_{\Phi,\omega,\alpha}^df(z))^p\omega(z)dV_\alpha(z)\right)^{1/p}\le C\left(\int_{\mathcal{H}}|f(z)|^p\omega(z)dV_\alpha(z)\right)^{1/p}.
\Ee
\elem
\begin{proof}
Using the last part in Lemma \ref{lem:covering} and that $\Phi\in B_p$, we obtain
\Beas
\|\mathcal{M}_{\Phi,\omega,\alpha}^df\|_{p,\omega,\alpha}^p &=& \int_0^\infty p\lambda^{p-1}\left|\{z\in \mathcal{H}:\mathcal{M}_{\Phi,\omega,\alpha}^df(z)>\lambda\}\right|_{\omega,\alpha}d\lambda\\ &\le& \int_0^\infty p\lambda^{p-1}\left(\int_{\{z\in \mathcal{H}:|f(z)|>\frac{\lambda}{2}\}}\Phi\left(\frac{|f(z)|}{\lambda}\right)\omega(z)dV_\alpha(z)\right)d\lambda\\ &=& p\int_{\mathcal{H}}\int_0^{2|f(z)|}\lambda^p\Phi\left(\frac{|f(z)|}{\lambda}\right)\frac{d\lambda}{\lambda}\omega(z)dV_\alpha(z)\\ &=& p\int_{\mathcal{H}}|f(z)|^p\left(\int_{1/2}^\infty \frac{\Phi(\lambda)}{\lambda^p}\frac{d\lambda}{\lambda}\right)\omega(z)dV_\alpha(z)\\ &=& C\int_{\mathcal{H}}|f(z)|^p\omega(z)dV_\alpha(z).
\Eeas
Here $C=p\int_{1/2}^\infty \frac{\Phi(\lambda)}{\lambda^p}\frac{d\lambda}{\lambda}$.
\end{proof}
Any weight in the B\'ekoll\`e-Bonami class, $\mathcal{B}_{p,\alpha}$, $1<p<\infty$ is an $\alpha$-doubling weight. 
\blem\label{lem:BekBonalphadoub}
Let $1<p<\infty$ and $\alpha>-1$. Then any weight $\omega\in \mathcal{B}_{p,\alpha}$ is $\alpha$-doubling, with doubling function $\varphi(t)=t^{p}$ and constant $K=[\omega]_{\mathcal{B}_{p,\alpha}}$.
\elem
\begin{proof}
Using H\"older's inequality and the definition of the class $\mathcal{B}_{p,\alpha}$, we obtain for any $E\subset Q_I$,
\Beas
\frac{|E|_\alpha}{|Q_I|_\alpha} &\le& \frac{\left(\int_E\omega(z)dV_\alpha(z)\right)^{1/p}\left(\int_{Q_I}\omega(z)^{-\frac{p'}{p}}dV_\alpha(z)\right)^{1/p'}}{|Q_I|_\alpha}\\ &\le& \left(\frac{|E|_{\omega,\alpha}}{|Q_I|_{\omega,\alpha}}\right)^{1/p}\frac{\left(\int_{Q_I}\omega(z)dV_\alpha(z)\right)^{1/p}\left(\int_{Q_I}\omega(z)^{-\frac{p'}{p}}dV_\alpha(z)\right)^{1/p'}}{|Q_I|_\alpha}\\ &\le& \left(\frac{|E|_{\omega,\alpha}}{|Q_I|_{\omega,\alpha}}[\omega]_{\mathcal{B}_{p,\alpha}}\right)^{1/p}.
\Eeas
Hence 
$$\frac{|Q_I|_{\omega,\alpha}}{|E|_{\omega,\alpha}}\le \left(\frac{|Q_I|_{\alpha}}{|E|_{\alpha}}\right)^p[\omega]_{\mathcal{B}_{p,\alpha}}.$$
\end{proof}

Let us now prove Theorem \ref{thm:main1}.
\begin{proof}[Proof of Theorem \ref{thm:main1}]
We start by observing that for any interval $I\subset \mathbb{R}$, $\|\chi_{Q_I}\|_{Q_{I},\Phi,\omega,\alpha}=1$. Thus taking $f=\chi_{Q_I}$, we obtain
$$\mu(Q_I)^{1/q}\le \left(\int_{Q_I}\left(\mathcal{M}_{\Phi,\omega,\alpha}(\chi_{Q_I})(z)\right)^qd\mu(z)\right)^{1/q}\le C\|\chi_{Q_I}\|_{p,\omega,\alpha}=C|Q_I|_{\omega,\alpha}^{1/p}$$
which leads to (\ref{eq:main12}). 

Conversely, assume that (\ref{eq:main12}) holds. Using Lemma \ref{lem:BekBonalphadoub}, Lemma \ref{lem:levelsetsembed}, we first obtain
\Beas
L &:=& \int_{\mathcal H}\left(\mathcal{M}_{\Phi,\omega,\alpha} f(z)\right)^qd\mu(z)\\ &=& \int_0^\infty q\lambda^{q-1}\mu\left(\{z\in \mathcal{H}:\mathcal{M}_{\Phi,\omega,\alpha}f(z)>\lambda\}\right)d\lambda\\ &\le& \int_0^\infty q\lambda^{q-1}\mu\left(\{z\in \mathcal{H}:\mathcal{M}_{\Phi,\omega,\alpha}^df(z)>\frac{\lambda}{C}\}\right)d\lambda\\ &\le& \int_0^\infty q\lambda^{q-1}\left|\{z\in \mathcal{H}:\mathcal{M}_{\Phi,\omega,\alpha}^df(z)>\frac{\lambda}{C}\}\right|_{\omega,\alpha}^{q/p}d\lambda\\ &\le& C^{q-p}\|\mathcal{M}_{\Phi,\alpha}^df\|_{p,\omega,\alpha}^{q-p}\int_0^\infty q\lambda^{p-1}\left|\{z\in \mathcal{H}:\mathcal{M}_{\Phi,\omega,\alpha}^df(z)>\frac{\lambda}{C}\}\right|_{\omega,\alpha}d\lambda\\ &=& C^q\frac{q}{p}\|\mathcal{M}_{\Phi,\omega,\alpha}^df\|_{p,\omega,\alpha}^q.
\Eeas
Where we have also used the inequality $$t^p \left|\{z\in \mathcal{H}:\mathcal{M}_{\Phi,\omega,\alpha}^df(z)>t\}\right|_{\omega,\alpha}\le \|\mathcal{M}_{\Phi,\alpha}^df\|_{p,\omega,\alpha}^p.$$

We easily conclude with the help of Lemma \ref{lem:estimmaxphi} that
\Beas \left(\int_{\mathcal H}\left(\mathcal{M}_{\Phi,\omega,\alpha} f(z)\right)^qd\mu(z)\right)^{1/q} &\le& C\|\mathcal{M}_{\Phi,\omega,\alpha}^df\|_{p,\omega,\alpha}\\ &\le& C\left(\int_{\mathcal{H}}|f(z)|^p\omega(z)dV_\alpha(z)\right)^{1/p}.
\Eeas
The proof is complete.
\end{proof}
In the case $1<p=q<\infty$, $\omega=1$ and $\mu=V_\alpha$, we obtain as a consequence, the following characterization of weighted Bergman spaces of the upper-half plane.
\bcor\label{cor:Bergequivdef}
Let $1<p<\infty$, and $\alpha>-1$. Assume $\Phi\in B_p$. Then for any analytic function $f$ on $\mathcal{H}$, the following are equivalent.
\begin{itemize}
\item[(i)] $f\in L^p(\mathcal{H}, dV_\alpha)$.
\item[(ii)] $\mathcal{M}_{\Phi,\alpha}f\in L^p(\mathcal{H}, dV_\alpha)$.
\end{itemize}
\ecor
\begin{proof}
That (i)$\Rightarrow$ (ii) is a special case of Theorem \ref{thm:main1}. To see that (ii)$\Rightarrow$ (i), observe that the Mean Value Theorem applied to the disc inscribed in the Carleson box $Q_I$ with centre $z\in \mathcal{H}$ allows one to see that there is a constant $C>0$ such that
$$|f(z)|\le \frac{C}{|Q_I|_\alpha}\int_{Q_I}|f(w)|dV_\alpha(w).$$ 
It is then enough to prove that there is a constant $K>0$ such that for any interval $I\subset \mathbb{R}$, $$\frac{1}{|Q_I|_\alpha}\int_{Q_I}|f(w)|dV_\alpha(w)\le K\|f\|_{Q_I,\Phi,\alpha}.$$
Using the convexity of $\Phi$, Jensen's inequality, and putting $\lambda=\|f\|_{Q_I,\Phi,\alpha}$, we easily obtain
\Beas
\Phi\left(\int_{Q_I}\frac{|f(w)|}{\lambda}\frac{dV_\alpha(w)}{|Q_I|_\alpha}\right) &\le& \int_{Q_I}\Phi\left(\frac{|f(w)|}{\lambda}\right)\frac{dV_\alpha(w)}{|Q_I|_\alpha}\le 1.
\Eeas 
Hence
$$\frac{1}{|Q_I|_\alpha}\int_{Q_I}|f(w)|dV_\alpha(w)\le \lambda\Phi(1)=\lambda.$$
The proof is complete.
\end{proof}
\subsection{A second proof of Theorem \ref{thm:main1}}
We consider the following system of dyadic grids,
$$\mathcal D^\beta:=\{2^j\left([0,1)+m+(-1)^j\beta\right):m\in \mathbb Z,\,\,\,j\in \mathbb Z \},\,\,\,\textrm{for}\,\,\,\beta\in \{0,1/3\}.$$
%For more on this system of dyadic grids and its applications, we refer to \cite{AlPottReg, HyPerez,Lerner,LerOmbroPerezetal,PR}. 
For $\beta=0$, $\mathcal D^0$ is the standard dyadic grid of $\mathbb R$, simply denoted $\mathcal D$.
\vskip .2cm 
We recall with \cite{PR} that given an interval $I\subset \mathbb R$, there is a dyadic interval $J\in \mathcal D^\beta$ for some $\beta\in \{0,1/3\}$ such that $I\subseteq J$ and $|J|\le 6|I|$. It follows that if $\omega$ is $\alpha$-doubling then in particular, we have that $|Q_J|_{\omega,\alpha}\le \varphi(6)K|Q_I|_{\omega,\alpha}$ where $\varphi$ and $K$ are given in (\ref{eq:alphadoubling}). Hence putting $\lambda:=\|f\|_{Q_J,\Phi,\omega,\alpha}$, we obtain 
\Beas\frac{1}{|Q_I|_{\omega,\alpha}}\int_{Q_I}\Phi\left(\frac{|f(z)|}{\varphi(6)K\lambda}\right)\omega(z)dV(z) &\le& \frac{\varphi(6)K}{|Q_J|_{\omega,\alpha}}\int_{Q_J}\Phi\left(\frac{|f(z)|}{\varphi(6)K\lambda}\right)\omega(z)dV(z)\\ &\le& \frac{1}{|Q_J|_{\omega,\alpha}}\int_{Q_J}\Phi\left(\frac{|f(z)|}{\lambda}\right)\omega(z)dV(z)\\ &\le& 1.
\Eeas
Thus  for any locally integrable function $f$,
\Be\label{eq:Maxfunctdyaineq}
\mathcal{M}_{\Phi,\omega,\alpha}f(z)\le C\sum_{\beta\in \{0,1/3\}}\mathcal{M}_{\Phi,\omega,\alpha}^{d,\beta} f(z),\,\,\,z\in \mathcal H
\Ee 
where $\mathcal{M}_{\Phi,\omega,\alpha}^{d,\beta}$ is defined as $\mathcal{M}_{\Phi,\omega,\alpha}$ but with the supremum taken only over dyadic intervals of the dyadic grid $\mathcal D^\beta$. When $\omega\equiv 1$, we use the notation $\mathcal{M}_{\Phi,\alpha}^{d,\beta}$, and if moreover, $\beta=0$, we just write $\mathcal{M}_{\Phi,\alpha}^{d}$. 
\vskip .2cm
It follows from the above observation that to prove that (\ref{eq:main11}) holds, it is enough to prove that this inequality holds for $\mathcal{M}_{\Phi,\omega,\alpha}^{d,\beta}$, $\beta=0,1/3$.

 We define $\mathcal{Q}^\beta$ by $$\mathcal{Q}^\beta:=\{Q=Q_I: I\in \mathcal{D}^\beta\}.$$
\begin{definition} Let $\alpha>-1$ and $\omega$ be a weight. For any $\gamma\geq 1$, a sequence of positive numbers $\{\lambda_Q\}_{Q\in \mathcal{Q}^\beta}$ is called a $(\omega,\alpha,\gamma)$-Carleson sequence, if there is a constant $C>0$ such that for any $R\in \mathcal{Q}^\beta$,
$$\sum_{Q\subseteq R,\,Q\in \mathcal{Q}^\beta}\lambda_Q\le C|R|_{\omega,\alpha}^\gamma.$$
\end{definition}

For $Q_I$, we denote by $T_I$ its upper half. That is $$T_I:=\{z=x+iy\in \mathbb C:x\in I\,\,\,\textrm{and}\,\,\,\frac{|I|}{2}<y<|I|\}.$$
%\vskip .3cm

Recall that for $1<p<\infty$ and $\alpha>-1$, a weight $\omega$ is said to belong to the B\'ekoll\`e-Bonami class $\mathcal{B}_{p,\alpha}$,  if $[\omega]_{\mathcal{B}_{p,\alpha}}<\infty$, where
$$[\omega]_{\mathcal{B}_{p,\alpha}}:=\sup_{I\subset \mathbb R,\,\,\, I\,\,\,\textrm{interval} }\left(\frac{1}{|Q_I|_\alpha}\int_{Q_I}\omega(z)dV_\alpha(z)\right)\left(\frac{1}{|Q_I|_\alpha}\int_{Q_I}\omega(z)^{1-p'}dV_\alpha(z)\right)^{p-1}.$$

For $p=\infty$, we say $\omega\in \mathcal{B}_{\infty,\alpha}$, if 
$$[\omega]_{\mathcal{B}_{\infty,\alpha}}:=\sup_{I\subset \mathbb R,\,\,\, I\,\,\,\textrm{interval} }\frac{1}{|Q_I|_{\omega,\alpha}}\int_{Q_I}\mathcal{M}_\alpha(\omega\chi_{Q_I})dV_\alpha(z)<\infty.$$
The following is proved in \cite{AlPottReg}.
\blem\label{lem:Binftyislarge}
Let $\alpha>-1$, and $1<p<\infty$. Then the class $\mathcal{B}_{\infty,\alpha}$ contains $\mathcal{B}_{p,\alpha}$. Moreover, for any weight $\omega\in \mathcal{B}_{p,\alpha}$,
$$[\omega]_{\mathcal{B}_{\infty,\alpha}}\le [\omega]_{\mathcal{B}_{p,\alpha}}.$$
\elem 
One can also prove that for any weight $\omega\in \mathcal{B}_{\infty,\alpha}$, the sequence $\{|Q_I|_{\omega,\alpha}\}_{I\in \mathcal{D}}$ is a $(\omega,\alpha,1)$-Carleson sequence.
\blem\label{lem:BinftyisCarl}
Let $\alpha>-1$, and  $\omega\in \mathcal{B}_{\infty,\alpha}$. Then for any $J\in \mathcal{D}^\beta$,
$$\sum_{I\subset J, I\in \mathcal{D}^\beta}|Q_I|_{\omega,\alpha}\le C_\alpha [\omega]_{\mathcal{B}_{\infty,\alpha}}|Q_I|_{\omega,\alpha}.$$
\elem
\begin{proof}
We have 
\Beas
\sum_{I\subset J, I\in \mathcal{D}^\beta}|Q_I|_{\omega,\alpha} &=& \sum_{I\subset J, I\in \mathcal{D}^\beta}\frac{|Q_I|_{\omega,\alpha}}{|Q_I|_\alpha}|Q_I|_\alpha\\ &=& C_\alpha\sum_{I\subset J, I\in \mathcal{D}^\beta}\frac{|Q_I|_{\omega,\alpha}}{|Q_I|_\alpha}|T_I|_\alpha\\ &\le& C_\alpha\sum_{I\subset J, I\in \mathcal{D}^\beta}\int_{T_I}\mathcal{M}_{\alpha}(\omega\chi_{Q_J})dV_\alpha\\ &\le& C_\alpha [\omega]_{\mathcal{B}_{\infty,\alpha}}|Q_I|_{\omega,\alpha}.
\Eeas
\end{proof}
\vskip .3cm
\btheo\label{thm:Carlembed0}
Let $\omega$ be a weight on $\mathcal{H}$ and $\gamma\ge 1$. Assume  $\{\lambda_{Q}\}_{Q\in \mathcal Q^\beta}$ is a sequence of positive numbers. Assume that
there exists some constant $A>0$ such that for any Carleson square $Q^0\in \mathcal Q^\beta$,
    \begin{equation*}\sum_{Q\subseteq Q^0, Q\in \mathcal Q^\beta}\lambda_{Q}\le A|Q^0|_{\omega,\alpha}^\gamma.\end{equation*}
Then there exists a constant $B>0$ such that for all $p\in (1,\infty)$, and any $\Phi\in B_p$,
\begin{equation*}\sum_{Q, Q\in \mathcal Q^\beta}\lambda_{Q}\|f\|_{Q,\Phi,\alpha}^{p\gamma}\le B \|\mathcal{M}_{\Phi,\omega, \alpha}^d\|_{p,\omega,\alpha}^{p\gamma}.\end{equation*}

\etheo

\begin{proof}
For simplicity of presentation, we assume that $\beta=0$.
We will also need the following inequality.
\Be\label{eq:weakineq}
\lambda^p|\{z\in \mathcal{H}:\mathcal{M}_{\Phi,\omega, \alpha}^df(z)>\lambda\}|_{\sigma,\alpha}\le \|\mathcal{M}_{\Phi,\omega, \alpha}^df\|_{p,\sigma,\alpha}^p.
\Ee

We can suppose that $f>0$. As in the case of $\Phi=1$ in \cite{HyPerez,Sehba}, we read $\sum_{Q\in \mathcal Q}\lambda_Q\|f\|_{Q,\Phi,\omega,\alpha}^{p\gamma}$ as an integral over the measure space $(\mathcal Q, \mu)$ built over the set of dyadic Carleson squares $\mathcal Q$, with $\mu$ the measure assigning to each square $Q\in \mathcal Q$ the measure $\lambda_Q$. Thus
\Beas
\sum_{Q\in \mathcal Q}\lambda_Q\|f\|_{Q,\Phi,\omega,\alpha}^{p\gamma} &=& \int_0^\infty p\gamma t^{p\gamma-1}\mu\left(\{Q\in \mathcal Q: \|f\|_{Q,\Phi,\omega,\alpha}>t\}\right)dt\\ &=& \int_0^\infty p\gamma t^{p\gamma-1}\mu(\mathcal Q_t)dt,
\Eeas
$\mathcal Q_t:=\{Q\in \mathcal Q: \|f\|_{Q,\Phi,\omega,\alpha}>t\}$. Let $\mathcal {Q}_t^*$ be the set of maximal dyadic Carleson squares $R$ with respect to the inclusion so that $\|f\|_{R,\Phi,\omega,\alpha}>t$. Then  $$\bigcup_{R\in
{\mathcal {Q}_t^*}}R=\{z\in \mathcal{H} :\mathcal{M}_{\Phi,\omega, \alpha}^df(z)>t\}.$$ It follows from the hypothesis on the sequence $\{\lambda_Q\}_{Q\in \mathcal Q}$ that
\Beas
\mu(\mathcal Q_t) &=& \sum_{Q\in \mathcal Q_t}\lambda_Q\le \sum_{R\in \mathcal {Q}_t^*}\sum_{Q\subseteq R}\lambda_Q\\ &\le& A\sum_{R\in \mathcal {Q}_t^*}|R|_{\omega,\alpha}^\gamma\le A\left(\sum_{R\in \mathcal {Q}_t^*}|R|_{\omega,\alpha}\right)^\gamma\\ &\le& A|\{z\in \mathcal{H} :\mathcal{M}_{\Phi,\omega, \alpha}^df(z)>t\}|_{\omega,\alpha}^\gamma.
\Eeas
Hence using (\ref{eq:weakineq}), we obtain
\Beas
S &:=& \sum_{Q\in \mathcal Q}\lambda_Q\|f\|_{Q,\Phi,\omega,\alpha}^{p\gamma}\\ &\le& A\int_0^\infty p\gamma t^{p\gamma-1}|\{z\in \mathcal{H} :\mathcal{M}_{\Phi,\omega, \alpha}^df(z)>t\}|_{\omega,\alpha}^\gamma dt\\ &=& A\int_0^\infty p\gamma t^{p-1}|\{z\in \mathcal{H} :\mathcal{M}_{\Phi,\omega, \alpha}^df(z)>t\}|_{\omega,\alpha}\\ && \left(t^p|\{z\in \mathcal{H} :\mathcal{M}_{\Phi,\omega, \alpha}^df(z)>t\}|_{\omega,\alpha}\right)^{\gamma-1} dt\\ &\le& A\gamma\|\mathcal{M}_{\Phi,\omega, \alpha}^d\|_{p,\omega,\alpha}^{p(\gamma-1)}\int_0^\infty p |\{z\in \mathcal{H} :\mathcal{M}_{\Phi,\omega, \alpha}^df(z)>t\}|_{\omega,\alpha}t^{p-1}dt\\ &\le& A\gamma\|\mathcal{M}_{\Phi,\omega, \alpha}^d\|_{p,\omega,\alpha}^{p\gamma}.
\Eeas
The proof is complete.
\end{proof}
We can now present our second proof of Theorem \ref{thm:main1}.
\begin{proof}[Proof of Theorem \ref{thm:main1}]
We only check that (\ref{eq:main11}) holds under (\ref{eq:main12}). As observed above, it is enough to prove that there exists a positive constant $C$ such that for any $\beta\in \{0, 1/3\}$, and any function $f$,
$$\int_{\mathcal H}(\mathcal{M}_{\Phi,\omega,\alpha}^{d,\beta}f(z))^{q}d\mu(z)\le C\left(\int_{\mathcal H}|f(z)|^{p}\omega(z)dV_\alpha(z)\right)^{q/p}.$$
This will follow from the following result.
\blem
Let $\alpha>-1$, $1<p\le q<\infty$, $\omega\in \mathbb{B}_{\infty.\alpha}$, and $\mu$ a positive Borel measure on $\mathcal H$. Assume that $\Phi\,\in B_{p}$ and that (\ref{eq:main12}) holds. Then there exists a positive constant $C$ such that for any $\beta\in \{0, 1/3\}$ and any function $f$,
$$\int_{\mathcal H}(\mathcal{M}_{\Phi,\omega,\alpha}^{d,\beta}f(z))^{q}d\mu(z)\le C\left(\int_{\mathcal H}|f(z)|^{p}\omega(z)dV_\alpha(z)\right)^{q/p}.$$
\elem
\begin{proof}
Let $a\ge 2$. To each integer $k$, we associate the level set
$$\Omega_{k}:=\{z\in \mathcal H: a^k<\mathcal{M}_{\Phi,\omega,\alpha}^{d,\beta}f(z)\leq a^{k+1}\}.$$
We observe with Lemma \ref{lem:covering} that 
$\Omega_{k}\subset \cup_{j=1}^{\infty}Q_{I_{k,j}},$ where
$Q_{I_{k,j}}$ ($I_{k,j}\in \mathcal{D}^{\beta}$) is a dyadic Carleson square maximal (with respect to the inclusion) such that
$$\|f\|_{Q_{I{k,j}},\Phi,\omega,\alpha}>a^k.$$

Using the above observations and condition (\ref{eq:main12}), we obtain  
 \Beas
 \int_{\mathcal H}(\mathcal{M}_{\Phi,\omega,\alpha}^{d,\beta}f(z))^{q}d\mu(z) &=& \sum_{k}\int_{\Omega_k}(\mathcal{M}_{\Phi,\omega,\alpha}^{d,\beta}f(z))^{q}d\mu(z)\\
 &\le& a^{q}\sum_{k}a^{kq}\mu(\Omega_k)\\
 &\le& a^{q}\sum_{k,j}a^{kq}\mu(Q_{I_{k,j}})\\
 &\le&
 a^{q}\sum_{k,j}\|f\|_{Q_{I{k,j}},\Phi,\omega,\alpha}^{q}\mu(Q_{I_{k,j}})\\ &\le& C\sum_{k,j}\|f\|_{Q_{I{k,j}},\Phi,\omega,\alpha}^{q}|Q_{I_{k,j}}|_{\omega,\alpha}^{q/p}\\ &\le& C\left(\sum_{k,j}\|f\|_{Q_{I{k,j}},\Phi,\omega,\alpha}^{p}|Q_{I_{k,j}}|_{\omega,\alpha}\right)^{q/p}.
\Eeas
As $\omega\in \mathbb{B}_{\infty,\alpha}$, we have that the sequence 
\Be\label{eq:pPhi} \lambda_Q= \left\{\begin{array}{lcr}|Q|_{\omega,\alpha} & \mbox{if} & Q=Q_{I_{k,j}}\,\,\,\textrm{for some}\,\,\,k,j\\ 0 & \mbox{otherwise} & \end{array}\right.\Ee
is a $(\omega,1,\alpha)$-Carleson sequence.
Hence using Theorem \ref{thm:Carlembed0} and Lemma \ref{lem:estimmaxphi}, we obtain

\Beas
 \int_{\mathcal H}(\mathcal{M}_{\Phi,\omega,\alpha}^{d,\beta}f(z))^{q}d\mu(z)  &\le& C\left(\sum_{k,j}\|f\|_{Q_{I{k,j}},\Phi,\omega,\alpha}^{p}|Q_{I_{k,j}}|_{\omega,\alpha}\right)^{q/p}\\ &\le& C\left(\int_{\mathcal H}(\mathcal{M}_{\Phi,\omega,\alpha}^{d,\beta}f(z))^{p}\omega(z)dV_\alpha(z)\right)^{q/p}\\ &\le& C\left(\int_{\mathcal H}|f(z)|^{p}\omega(z)dV_\alpha(z)\right)^{q/p}.
\Eeas
\end{proof}
The proof is complete.
\end{proof}
\section{Proof of Theorem \ref{thm:main2}}

Let us start by proving the following result.
%$(\textrm{ii})\Rightarrow (\textrm{i})$, it is enough by the observations made at the beginning of this section to prove the following.
\blem\label{lem:main11}
Let $1< q<p<\infty$. Let $\omega\in \mathbb{B}_{\infty,\alpha}$, and $\mu$ be a positive measure on $\mathcal{H}$. Assume that $\Phi\,\in B_{p}$ and that the function
\Be\label{eq:main111} K_{\mu}(z):=\sup_{I\subset \mathbb R,\,\,\,z\in Q_I}\frac{\mu(Q_{I})}{|Q_I|_{\omega,\alpha}}
\Ee belongs to
$L^s(\mathcal H,\omega dV_\alpha)$ where $s=\frac{p}{p-q}$.  Then there is a positive constant $C>0$ such that for any $f\in L^p(\mathcal H, \omega dV_\alpha)$, and any $\beta\in \{0,\frac{1}{3}\}$,
\Be\label{eq:main112}
\left(\int_{\mathcal H}(\mathcal{M}_{\Phi,\omega,\alpha}^{d,\beta} f(z))^qd\mu(z)\right)^{1/q}\le C\left(\int_{\mathcal H}|f(z)|^p\omega(z)dV_\alpha(z)\right)^{1/p}.
\Ee
\elem
\begin{proof}
Let $a\ge 2$. To each integer $k$, we associate the set
$$\Omega_{k}:=\{z\in \mathcal H: a^k<\mathcal{M}_{\Phi,\omega,\alpha}^{d,\beta}f(z)\leq a^{k+1}\}.$$
We recall with Lemma \ref{lem:covering} that 
$\Omega_{k}\subset \cup_{j=1}^{\infty}Q_{I_{k,j}},$ where
$Q_{I_{k,j}}$ ($I_{k,j}\in \mathcal{D}^{\beta}$) is a dyadic cube maximal (with respect to the inclusion) such that
$$\|f\|_{Q_{I{k,j}},\Phi,\omega,\alpha}>a^k.$$

For $\beta\in \{0,\frac{1}{3}\}$, we define
$$K_{d,\mu}^\beta(z):=\sup_{I\in \mathcal{D}^\beta,z\in Q_I}\frac{\mu(Q_{I})}{\omega(Q_{I})}.$$ Then $K_{d,\mu}^\beta(z)\leq K_\mu(z)$ for any $z\in \mathcal{H}$. Hence, that $K_\mu(z)$ belongs to $L^{p/(p-q)}(\mathcal H,\omega dV_\alpha)$ implies that $K_{d,\mu}^\beta\in L^{p/(p-q)}(\mathcal H,\omega dV_\alpha)$.
\vskip .3cm
Proceeding as in the second proof of Theorem \ref{thm:main1}, we obtain at first that 
 \Beas
 \int_{\mathcal H}(\mathcal{M}_{\Phi,\omega,\alpha}^{d,\beta}f(z))^{q}d\mu(z) &=& \sum_{k}\int_{\Omega_k}(\mathcal{M}_{\Phi,\omega,\alpha}^{d,\beta}f(z))^{q}d\mu(z)\\
  &\le&
 a^{q}\sum_{k,j}\|f\|_{Q_{I{k,j}},\Phi,\omega,\alpha}^{q}\mu(Q_{I_{k,j}}).
\Eeas
Now using H\"older's inequality, we obtain 
\Beas
\int_{\mathcal H}(\mathcal{M}_{\Phi,\omega,\alpha}^{d,\beta}f(z))^{q}d\mu(z) &\le& a^{q}\sum_{k,j}\|f\|_{Q_{I{k,j}},\Phi,\omega,\alpha}^{q}\mu(Q_{I_{k,j}})\\
 &=& a^{q}\sum_{k,j}\|f\|_{Q_{I{k,j}},\Phi,\omega,\alpha}^{q}
 \frac{\mu(Q_{I_{k,j}})}{|Q_{I_{k,j}}|_{\omega,\alpha}}|Q_{I_{k,j}}|_{\omega,\alpha}\\
&\le& a^qA^{q/p}B^{1/s}
\Eeas
where $$A:=\sum_{k,j}\|f\|_{Q_{I{k,j}},\Phi,\omega,\alpha}^{p}|Q_{I_{k,j}}|_{\omega,\alpha}$$
and
$$B:=\sum_{k,j}\left(\frac{\mu(Q_{I_{k,j}})}{|Q_{I_{k,j}}|_{\omega,\alpha}}\right)^{s}|Q_{I_{k,j}}|_{\omega,\alpha}.$$
The estimate of $A$ is already obtained in the second proof of Theorem \ref{thm:main1}.
Let us estimate $B$. As $\omega$ is $\alpha$-doubling, we obtain
\Beas
B &:=&\left(\sum_{k,j}\left(\frac{\mu(Q_{I_{k,j}})}{|Q_{I_{k,j}}|_{\omega,\alpha}}\right)^{s}|Q_{I_{k,j}}|_{\omega,\alpha}\right)\\
%\left(\sum_{k,j}\left(\frac{\mu(Q_{I_{k,j}})}{\omega(Q_{I_{k,j}})}\right)^{s}\omega(Q_{I_{k,j}})\right)^{\frac{1}{s}}&=\nonumber
&\lesssim&
\sum_{k,j}\left(\frac{\mu(Q_{I_{k,j}})}{|Q_{I_{k,j}}|_{\omega,\alpha}}\right)^{s}|T_{I_{k,j}}|_{\omega,\alpha}\\
&\lesssim& \sum_{k,j}\int_{T_{I_{k,j}}}\left(\frac{\mu(Q_{I_{k,j}})}{|Q_{I_{k,j}}|_{\omega,\alpha}}\right)^{s}\omega(z)dV_\alpha(z)\\
&\lesssim& \sum_{k,j}\int_{T_{I_{k,j}}}\left(K_{d,\mu}^\beta(z)\right)^{s}\omega(z)dV_\alpha(z)\\
&\lesssim&
\int_{\mathcal H}\left(K_{d,\mu}^\beta(z)\right)^{s}\omega(z)dV_\alpha(z)=C\|K_{d,\mu}^\beta\|_{s,\omega,\alpha}^s<\infty.
\Eeas

The proof of the lemma is complete.
\end{proof}
We can now prove the Theorem \ref{thm:main2}.
\begin{proof}[Proof of Theorem~{\rm\ref{thm:main2}}]
The proof of the sufficiency follows from Lemma \ref{lem:main11} and the observations made at the beginning of this section. Let us prove the necessity. For this, we first check that the Hardy-Littlewood maximal function is pointwise dominated by the generalized maximal function. Indeed, let $I\subset \mathbb{R}$ be an interval an put 
$\lambda:=\|f\|_{Q_I,\Phi,\omega,\alpha}$. Using the convexity of $\Phi$, we obtain 
\Beas\Phi\left(\frac{1}{|Q_I|_{\omega,\alpha}}\int_{Q_I}\frac{|f(z)|}{\lambda}\omega(z)dV(z)\right) &\le& \frac{1}{|Q_I|_{\omega,\alpha}}\int_{Q_J}\Phi\left(\frac{|f(z)|}{\lambda}\right)\omega(z)dV(z)\\  &\le& 1.
\Eeas
Thus for any $z\in \mathcal{H}$, and for any interval $I\subset \mathbb{R}$ such that $z\in Q_I$,
$$\frac{1}{|Q_I|_{\omega,\alpha}}\int_{Q_I}|f(w)|\omega(w)dV(w)\le \|f\|_{Q_I,\Phi,\omega,\alpha}.$$
Hence $$\mathcal{M}_{\omega,\alpha}f(z)\le \mathcal{M}_{\Phi,\omega,\alpha}f(z)$$
for any locally integrable function $f$. It follows that if (\ref{eq:main11}) holds for the operator $\mathcal{M}_{\Phi,\omega,\alpha}$, it also holds for $\mathcal{M}_{\Phi,\alpha}$. That (\ref{eq:main11}) holds for $\mathcal{M}_{\Phi,\alpha}$ implies that the function given by (\ref{eq:main21}) belongs to $L^{p/(p-q)}(\mathcal{H}, \omega dV_\alpha)$ is proved in \cite{CarnotBenoit}.
The proof is complete.
\end{proof}
\section{Proof of Theorem \ref{thm:main3}}
We will be using discretization once more.
\begin{proof}
As seen before, we only need to establish the inequality (\ref{eq:main32}) for the dyadic maximal function $\mathcal{M}_\alpha^{d,\beta}$, $\beta\in \{0, 1/3\}$.
\vskip .3cm
Let $a\ge 2$. To each integer $k$, we associate the set
$$\Omega_{k}:=\{z\in \mathcal H: a^k<\mathcal{M}_{\alpha}^{d,\beta}f(z)\leq a^{k+1}\}.$$
As special case of Lemma \ref{lem:covering}, we have that
$\Omega_{k}\subset \cup_{j=1}^{\infty}Q_{I_{k,j}},$ where
$Q_{I_{k,j}}$ ($I_{k,j}\in \mathcal{D}^{\beta}$) is a dyadic cube maximal (with respect to the inclusion) such that
$$\frac{1}{|Q_{I_{k,j}}|_\alpha}\int_{Q_{I_{k,j}}}|f(z)|dV_\alpha(z)>a^k.$$
Proceeding as in the second proof of Theorem \ref{thm:main1}, we obtain
\Beas
 \int_{\mathcal H}(\mathcal{M}_{\alpha}^{d,\beta}f(z))^{q}d\mu(z) &=& \sum_{k}\int_{\Omega_k}(\mathcal{M}_{\alpha}^{d,\beta}f(z))^{q}d\mu(z)\\
 &\le& a^{q}\sum_{k}a^{kq}\mu(\Omega_k)\\
 &\le& a^{q}\sum_{k,j}a^{kq}\mu(Q_{I_{k,j}})\\
 &\le&
 a^{q}\sum_{k,j}\left(\frac{1}{|Q_{I_{k,j}}|_\alpha}\int_{Q_{I_{k,j}}}|f(z)|dV_\alpha(z)\right)^q\mu(Q_{I_{k,j}}).
\Eeas
Now let $\Psi$ be the complementary function of the Young function $\Phi$. Recall the following H\"olders's inequality:
\Be\label{eq:holderhgene}
\frac{1}{|Q_{I}|_\alpha}\int_{Q_{I}}|(fg)(z)|dV_\alpha(z)\le \|f\|_{Q_I,\Phi,\alpha}\|g\|_{Q_I,\Psi,\alpha}
\Ee
 Using the above generalized H\"older's inequality and (\ref{eq:main31}), we obtain
\Beas
 \int_{\mathcal H}(\mathcal{M}_{\alpha}^{d,\beta}f(z))^{q}d\mu(z) 
 &\le&
 a^{q}\sum_{k,j}\left(\frac{1}{|Q_{I_{k,j}}|_\alpha}\int_{Q_{I_{k,j}}}|f(z)|dV_\alpha(z)\right)^q\mu(Q_{I_{k,j}})\\ &=& a^{q}\sum_{k,j}\left(\frac{1}{|Q_{I_{k,j}}|_\alpha}\int_{Q_{I_{k,j}}}|f(z)|\omega(z)\omega(z)^{-1}dV_\alpha(z)\right)^q\mu(Q_{I_{k,j}})\\ &\le& a^{q}\sum_{k,j}\|f\omega\|_{Q_{I_{k,j}},\Phi,\alpha}^q\|\omega^{-1}\|_{Q_{I_{k,j}},\Psi,\alpha}^q\mu(Q_{I_{k,j}})\\ &\le& C\sum_{k,j}\|f\omega\|_{Q_{I_{k,j}},\Phi,\alpha}^q|Q_{I_{k,j}}|_\alpha^{q/p}.
\Eeas
It follows using Lemma \ref{lem:estimmaxphi} that
\Beas
 \int_{\mathcal H}(\mathcal{M}_{\alpha}^{d,\beta}f(z))^{q}d\mu(z) 
 &\le&
 C\sum_{k,j}\|f\omega\|_{Q_{I_{k,j}},\Phi,\alpha}^q|Q_{I_{k,j}}|_\alpha^{q/p}\\ &\le& C\left(\sum_{k,j}\|f\omega\|_{Q_{I_{k,j}},\Phi,\alpha}^p|Q_{I_{k,j}}|_\alpha\right)^{q/p}\\ &\le& C\left(\sum_{k,j}\|f\omega\|_{Q_{I_{k,j}},\Phi,\alpha}^p|T_{I_{k,j}}|_\alpha\right)^{q/p}\\ &=& C\left(\sum_{k,j}\int_{T_{I_{k,j}}}\|f\omega\|_{Q_{I_{k,j}},\Phi,\alpha}^pdV_\alpha\right)^{q/p}\\ &\le& C\left(\int_{\mathcal{H}}(\mathcal{M}_{\Phi,\alpha}^{d,\beta}(\omega f)(z))^{p}dV_\alpha\right)^{q/p}\\ &\le& C\left(\int_{\mathcal{H}}|(\omega f)(z)|^{p}dV_\alpha(z)\right)^{q/p}.
\Eeas
The proof is complete.
\end{proof}
\section{Proof of Theorem \ref{thm:main4} and Theorem \ref{thm:main5}}
To prove the inequality (\ref{eq:main41}), again, we only need to prove that the same inequality holds when $\mathcal{M}_{\Phi,\alpha}$ is replaced by its dyadic counterparts $\mathcal{M}_{\Phi,\alpha}^{d,\beta}$, $\beta=0,\frac{1}{3}$. Recall that if $a\ge 2$, then we associate to each integer $k$,  the set
$$\Omega_{k}:=\{z\in \mathcal H: a^k<\mathcal{M}_{\Phi,\alpha}^{d,\beta}f(z)\leq a^{k+1}\}$$
and we have that
$\Omega_{k}\subset \cup_{j=1}^{\infty}Q_{I_{k,j}},$ where
$Q_{I_{k,j}}$ ($I_{k,j}\in \mathcal{D}^{\beta}$) is a dyadic cube maximal (with respect to the inclusion) such that
$$\|f\|_{Q_{I_{k,j}},\Phi,\alpha}>a^k.$$ 
\begin{proof}[Proof of Theorem \ref{thm:main4}]
Following the same decomposition as in the proof of Theorem \ref{thm:main1} and using the assumption on the measure $\mu$, we obtain
\Beas
 \int_{\mathcal H}(\mathcal{M}_{\Phi,\alpha}^{d,\beta}f(z))^{q}d\mu(z) &=& \sum_{k}\int_{\Omega_k}(\mathcal{M}_{\Phi,\alpha}^{d,\beta}f(z))^{q}d\mu(z)\\
  &\le&
 C\sum_{k,j}\|f\|_{Q_{I_{k,j}},\Phi,\alpha}^q\mu(Q_{I_{k,j}})\\ &\le& C\sum_{k,j}\|f\|_{Q_{I_{k,j}},\Phi,\alpha}^q|Q_{I_{k,j}}|_{\omega,\alpha}^{q/p}\\ &=& C\sum_{k,j}\|f\|_{Q_{I_{k,j}},\Phi,\alpha}^q\left(\frac{|Q_{I_{k,j}}|_{\omega,\alpha}}{|Q_{I_{k,j}}|_\alpha}\right)^{q/p}|Q_{I_{k,j}}|_\alpha^{q/p}\\ &=& C\sum_{k,j}\|f\left(\frac{|Q_{I_{k,j}}|_{\omega,\alpha}}{|Q_{I_{k,j}}|_\alpha}\right)^{1/p}\|_{Q_{I_{k,j}},\Phi,\alpha}^q|Q_{I_{k,j}}|_\alpha^{q/p}.
\Eeas
This leads us to
\Beas
 \int_{\mathcal H}(\mathcal{M}_{\Phi,\alpha}^{d,\beta}f(z))^{q}d\mu(z)  &\le& C\left(\sum_{k,j}\|f\left(\frac{|Q_{I_{k,j}}|_{\omega,\alpha}}{|Q_{I_{k,j}}|_\alpha}\right)^{1/p}\|_{Q_{I_{k,j}},\Phi,\alpha}^p|Q_{I_{k,j}}|_\alpha\right)^{q/p}\\ &\le& C\left(\sum_{k,j}\|f\left(\frac{|Q_{I_{k,j}}|_{\omega,\alpha}}{|Q_{I_{k,j}}|_\alpha}\right)^{1/p}\|_{Q_{I_{k,j}},\Phi,\alpha}^p|T_{I_{k,j}}|_\alpha\right)^{q/p}\\ &\le& C\left(\int_{\mathcal H}\left(\mathcal{M}_{\Phi,\alpha}^{d,\beta}\left(f(\mathcal{M}_\alpha\omega)^{1/p}\right)(z)\right)^{p}dV_{\alpha}(z)\right)^{q/p}\\ &\le& C\left(\int_{\mathcal H}|f(z)|^p\mathcal{M}_\alpha\omega(z)dV_{\alpha}(z)\right)^{q/p}.
\Eeas
\end{proof}
\begin{proof}[Proof of Theorem \ref{thm:main5}]
Using the same notations as above, we first obtain
\Beas
 \int_{\mathcal H}(\mathcal{M}_{\Phi,\alpha}^{d,\beta}f(z))^{q}d\mu(z) &=& \sum_{k}\int_{\Omega_k}(\mathcal{M}_{\Phi,\alpha}^{d,\beta}f(z))^{q}d\mu(z)\\
  &\le&
 C\sum_{k,j}\|f\|_{Q_{I_{k,j}},\Phi,\alpha}^q\mu(Q_{I_{k,j}})\\ &=& C\sum_{k,j}\|f\|_{Q_{I_{k,j}},\Phi,\alpha}^q\frac{\mu(Q_{I_{k,j}})}{|Q_{I_{k,j}}|_{\omega,\alpha}}|Q_{I_{k,j}}|_{\omega,\alpha}.
 \Eeas
An easy application of H\"older's inequality to the last term in right of the above inequalities gives us 
\Beas
 \int_{\mathcal H}(\mathcal{M}_{\Phi,\alpha}^{d,\beta}f(z))^{q}d\mu(z) &\le& C\sum_{k,j}\|f\|_{Q_{I_{k,j}},\Phi,\alpha}^q\frac{\mu(Q_{I_{k,j}})}{|Q_{I_{k,j}}|_{\omega,\alpha}}|Q_{I_{k,j}}|_{\omega,\alpha}\\ &\le& CA^{q/p}B^{\frac{p-q}{p}}.
 \Eeas
 where 
 $$A:=\sum_{k,j}\|f\|_{Q_{I_{k,j}},\Phi,\alpha}^p|Q_{I_{k,j}}|_{\omega,\alpha}$$
 and 
 $$B:=\sum_{k,j}\left(\frac{\mu(Q_{I_{k,j}})}{|Q_{I_{k,j}}|_{\omega,\alpha}}\right)^{\frac{p}{p-q}}|Q_{I_{k,j}}|_{\omega,\alpha}.$$
 Following the lines of the proof of Theorem \ref{thm:main4}, we see that
 $$A\le C\int_{\mathcal H}|f(z)|^p\mathcal{M}_\alpha\omega(z)dV_{\alpha}(z).$$
The estimate of the term $B$ is quite harmless, we easily obtain
\Beas
B &:=& \sum_{k,j}\left(\frac{\mu(Q_{I_{k,j}})}{|Q_{I_{k,j}}|_{\omega,\alpha}}\right)^{\frac{p}{p-q}}|Q_{I_{k,j}}|_{\omega,\alpha}\\ &=& \sum_{k,j}\left(\frac{\mu(Q_{I_{k,j}})}{|Q_{I_{k,j}}|_{\omega,\alpha}}\right)^{\frac{p}{p-q}}\frac{|Q_{I_{k,j}}|_{\omega,\alpha}}{|Q_{I_{k,j}}|_\alpha}|Q_{I_{k,j}}|_\alpha\\ &\le & \sum_{k,j}\left(\frac{\mu(Q_{I_{k,j}})}{|Q_{I_{k,j}}|_{\omega,\alpha}}\right)^{\frac{p}{p-q}}\frac{|Q_{I_{k,j}}|_{\omega,\alpha}}{|Q_{I_{k,j}}|_\alpha}|T_{I_{k,j}}|_\alpha\\ &\le& \sum_{k,j}\int_{T_{I_{k,j}}}\left(\frac{\mu(Q_{I_{k,j}})}{|Q_{I_{k,j}}|_{\omega,\alpha}}\right)^{\frac{p}{p-q}}\frac{|Q_{I_{k,j}}|_{\omega,\alpha}}{|Q_{I_{k,j}}|_\alpha}dV_\alpha(z)\\ &\le& \int_{\mathcal{H}} (K_{\mu}(z))^{\frac{p}{p-q}}\mathcal{M}_\alpha\omega(z)dV_\alpha(z).
\Eeas  
The proof is complete.
\end{proof}

\end{document}